\documentclass{pspum-l}
\usepackage{
  amsmath,epsfig,amscd,xspace,stmaryrd}
\usepackage{enumitem}
\usepackage[usenames]{color}
\usepackage[english]{babel}                   
\usepackage[latin1]{inputenc} 
\usepackage{amsmath,amssymb}
\usepackage{epsfig}
\usepackage[all]{xy}
\xyoption{2cell}
\usepackage{theorem}
\usepackage{framed}

\newdir{ )}{{}*!/-5pt/@^{(}}
\newcommand{\eg}{e.g.,\xspace}
\newcommand{\ie}{i.e.,\xspace}

\CompileMatrices

\definecolor{DarkBlue}{rgb}{0,0.06,0.25}

\newcommand{\Kinv}{\text{inv}}

\newcommand{\Z}{{\mathbb Z}}
\newcommand{\Q}{{\mathbb Q}}

\newcommand{\C}{{\mathbb C}}
\newcommand{\F}{{\mathbb F}}

\newcommand{\THH}{\mathrm{T\hspace{-.5mm}H\hspace{-.5mm}H}}
\newcommand{\HH}{\mathrm{H\hspace{-.5mm}H}}
\newcommand{\HC}{\mathrm{H\hspace{-.5mm}C}}
\newcommand{\TC}{\mathrm{T\hspace{-.5mm}C}}
\newcommand{\TR}{\mathrm{T\hspace{-.5mm}R}}

\newcommand{\fib}{\twoheadrightarrow}

\newcommand{\we}{\overset{\sim}{\rightarrow}}

\newcommand{\smsh}{{\wedge}}

\newcommand{\ess}{\mathbb S}

\newcommand{\T}{\mathbb T} 
\theoremstyle{plain}

\def\qedbox{$\Box$}
\newbox\QedBox \setbox\QedBox\vbox{\hrule height 1ex width .618ex}
\def\qedbox{\copy\QedBox}

\makeatletter
\def\qed{%
   {%
      \unskip
      \nobreak \hfil
      \penalty 50               
      \hskip 3em                
      \null \nobreak \hfil
      \qedbox
      \parfillskip=\z@skip
      \finalhyphendemerits=\z@  
      \endgraf                  
   }}
\makeatother

\usepackage{hyperref}
\begin{document}
\title{Applications of topological cyclic homology to algebraic K-theory}
\author{Bj\o rn Ian Dundas}

\maketitle
\begin{quote}{\small
    {\textbf Abstract.}
    Algebraic K-theory has applications in a broad range of mathematical subjects, from number theory to functional analysis.
    It is also fiendishly hard to calculate.
    Presently there are two main inroads: motivic homotopy theory and cyclic homology.
    I've been asked to present an overview of the applications of topological cyclic homology to algebraic K-theory ``from a historical perspective''.
    The timeline spans from the very early days of algebraic K-theory to the present, starting with ideas in the seventies around the ``tangent space'' of algebraic K-theory all the way to the current state of affair where we see a resurgence in structural theorems, calculations and a realization that variants of cyclic homology have important things to say beyond the moorings to K-theory.
}
\end{quote}

In what follows I try to give a quasi-historical overview of topological cyclic homology, $\TC$, from the perspective of algebraic K-theory by selecting some papers I think are good illustrations.  These papers may or may not be the ``first'' or the ``best'' in any sense and since I've chosen them there is some bias, but regardless, by starting from these and see who cites them and who they cite I hope it should be possible to get a fairly broad grasp of the subject.
I have largely ignored development that is  not about the {\bf connection} between the two theories.
To cite C.~Weibel \cite{MR1732049} ``Like all papers, this one has a finiteness obstruction [\dots]
We apologize in advance to the large number of people whose important contributions we have skipped over.''

Before we start on the historical development, there are some ICM/ECM addresses and books that should be mentioned for reference: M.~Karoubi \cite{MR913964} (which was my main source as a student),
T.~G.~Goodwillie \cite{MR1159249}, I.~Madsen \cite{MR1474979}, \cite{MR1341845},
J.-L.~Loday \cite{MR1600246} (the readers who wish to check the definitions of more algebraic matters like K\"ahler differentials, Hochschild homology and so on should have this handy since it would break the flow to do so here), L.~Hesselholt \cite{MR1957052},  B.~I.~Dundas, Goodwillie and R.~McCarthy \cite{MR3013261} and J.~Rognes \cite{MR3728661}.

Copyright conundrums rob this paper for the most enjoyable part of the original talk [from which this paper was much too hastily assembled]:
the happy faces of many of my heroes and fellow mathematicians.  Please picture them yourself or do a web search for the talk.

I'll use the following shorthand: K-theory means algebraic K-theory, ring spectrum and $\ess$-algebra means $E_1$-ring spectrum and commutative ring spectrum means $E_\infty$-ring spectrum

I thank the referee for many good suggestions.  At certain points the suggestions have forced me to abandon the idea of largely leaving out material that had not yet  made it to publication by the time I gave the talk.  This risks being partial to the emerging papers that I know about at the expense of some potentially breathtaking ones that have not yet caught my attention.  For this I apologize.

\section{Prehistory}
\label{sec:prehistory}
The history of homological invariants of 
K-theory 
goes far back, but for our purposes a natural place to start is with S.~M.~Gersten's $\mathrm d\log$ map
and S.~Bloch's interpretation \cite{MR0466264} of this map as providing information about the ``tangent space'' of algebraic K-theory. This can be seen as a beginning of the story for many reasons and it actually appears in the 1972 Battelle  proceedings in which D.~Quillen first develops fully his algebraic K-theory.

If $A$ is a commutative ring and $\Omega^*_A$ its ring of K\"ahler differentials, Gersten defines a graded ring homomorphism
$$\mathrm d\log\colon K_*(A)\to\Omega^*_A$$
via Grothendieck's theory of Chern classes,
extending the map $\mathrm{GL}_1(A)\to \Omega^1_A$ sending a unit $g$ to $\frac1g\,dg$.  Bloch defines the \emph{tangent space} $TK_n(A)$ of K-theory at $A$ to be the kernel of the map $K_n(A[\epsilon]/\epsilon^2)\to K_n(A)$ 
induced by 
$\epsilon\mapsto0$ and proves that $\Omega^{n-1}_A$ splits off $TK_n(A)$ if $A$ is a local ring with $2$ invertible.  When $n=2$ this was known to van~der~Kallen, and one can think of this as a first step towards identifying relative K-theory in terms of cyclic homology-like invariants.

Another observation in the seventies with implications to $\TC$ is K.~Dennis' trace map to Hochschild homology
$$K_n(A)\to \HH_n(A),$$
related to the Hattori-Stallings trace \cite{MR175950}, \cite{MR202807}.  Simultaneously, there were deep investigations by many authors of the connections to crystalline cohomology, Chow groups, algebraic cycles, Chern classes and so on.

A quick word about how the maps of Gersten and Dennis fit in the picture to emerge may be in order.  Thinking of K-theory as having something to do with vector bundles and K\"ahler differentials with differential forms, Gersten's map
seems very geometric.  Upon staring at a concrete formula, the  Dennis trace map may seem more algebraic.  This is an illusory distinction.  In the smooth case the Hochschild-Kostant-Rosenberg theorem asserts that the ``antisymmetrization map'' $\epsilon\colon\Omega^n_A\to\HH_n(A)$ is an isomorphism, and -- getting a bit ahead of ourselves -- for general commutative $A$ the de~Rham differential $d\colon\Omega^n_A\to\Omega^{n+1}_A$ corresponds under $\epsilon$ to Connes' operator $B\colon\HH_n(A)\to\HH_{n+1}(A)$ coming from the circle action on $\HH$ which will be central when we get to cyclic homology. Connes' idea is that with Hochschild homology we gain access also to ``non-commutative geometry''.   For more background and definitions the reader may consult \cite{MR1600246} and \cite{MR1303779}.

In the second half of the seventies F.~Waldhausen enters the fray in full force.  He is motivated by the study of manifolds and the prospect that algebraic K-theory may have much more to offer in this direction if one allows a generalization of rings he envisions as ``Brave New Rings''.  In modern parlance, he is considering algebras over the sphere spectrum $\ess$, and in particular spherical group rings $\ess[G]$ whose algebraic K-theory gives deep information about the classifying space $BG$ of the (simplicial or topological) group $G$.
Through work of for instance J.~H.~C.~Whitehead and C.~T.~C.~Wall and the famous s-cobordism theorem of Barder-Mazur-Stallings it was known that important information about a pointed connected space $X$ could be obtained by looking at the K-theory of the integral group ring $\Z[\pi_1X]$.  However, it was clear that much information was lost.
It was primarily the lower K-groups that carried valuable information and even that picture was marred by fringe effects and missing torsion information.  Replacing the fundamental group by the loop group is a good idea, but Waldhausen's daring leap is to realize that the relevant questions involve not linear algebra over the integers, but in a stable context.  Hence the focus shifts to spherical group rings.  The interested reader may consult the book by Waldhausen, Jahren and Rognes \cite{MR3202834} and its references for the geometric consequences, in particular \cite{MR932266} and \cite{MR921486}.
We will have a few more details to offer when we revisit this theme in Section~\ref{sec:rat} (rational results) and Section~\ref{sec:grouprings} (integral results) where the relation with cyclic homology is at the forefront.



In \cite{MR561230} Waldhausen describes how the Dennis trace for a ring $A$ factors through a universal homology theory under K-theory he dubs \emph{stable K-theory}
$$K(A)\to K^S(A)\to\HH(A)$$
(here and later $K(A)$ is a naturally chosen spectrum with $\pi_nK(A)=K_n(A)$ -- this turns out to be an important shift of paradigm: the groups are no longer our main target).
The first map makes sense also for $A$ an $\ess$-algebra, and for spherical group rings it splits (with a geometrically important complement), 
and a natural question is what $\HH$ should mean for $\ess$-algebras.
Eventually Waldhausen and Goodwillie conclude that there ought to exist a Hochschild homology over the sphere spectrum with tensor replaced by the (at the time still conjectural) smash product of spectra.  They call this theory \emph{topological Hochschild homology} $\THH$ and conjecture that stable K-theory and topological Hochschild homology are equivalent.

\section{Cyclic homology in the 80's}
\label{sec:cyclichomology}

In the early eighties A. Connes defined \emph{cyclic homology}, $\HC$ \cite{ConnesExt}, \cite{MR823176}, a de~Rham-like theory for associative rings which can also be seen as the homotopy orbits $\HH(A)_{h\T}$ by the cyclic action on the standard representation of Hochschild homology, see \cite{MR1600246}.  Although not Connes' main focus at the time, to us the central questions at this juncture might be
\begin{itemize}[itemsep=0pt]
\item what is the relation between algebraic K-theory and cyclic homology?
\item in particular, is there a ``$\mathrm d\log$'' like Gersten's or more generally some sort of Chern theory for associative rings and how does it compare with Bloch's ``tangent space'' insight?
\item How does cyclic homology fit in Waldhausen's program?  Is there a version for $\ess$-algebras and how much information about K-theory can be retained?
\end{itemize}
The first years of cyclic homology see an explosion of ideas appearing simultaneously.  Independently, B. Tsygan  \cite{MR695483} had discovered the theory and dubbed it \emph{additive K-theory} since (at least rationally) it is ``like K-theory but with the Lie algebra $\frak{gl}(A)$ playing the r\^ole of the general linear group $\mathrm{GL}(A)$''; a point of view investigated both in Feigin-Tsygan \cite{MR923136} and Loday-Quillen \cite{MR780077}.  This offers a tantalizing geometric perspective: if we pretend that $\mathrm{GL}_n(A)$ is a Lie group $G$, then the Lie group $\frak{gl}_n(A)$ should be the tangent space $T_1G$ at the unit $1\in G$ and the inverse of the exponential map
$T_1G\to G$
should be a local isomorphism connecting K-theory (made from $\mathrm{GL_n}$) to cyclic homology (made from $\frak{gl}_n$).  A na\"\i ve look at the formula $\log(1-x)=\sum_{j>0}\frac{x^j}j$ gives the impression that this idea only works over the rationals (we need to divide out by the denominators) and somehow $(x^j)_{j>0}$ should converge rather rapidly to zero.  This fits hand in glove with Bloch's approach: he considered the situation $x^2=0$ (to interpret the tangent space).


\subsection{Rational results}
\label{sec:rat}
In hindsight it is clear that the some of the ideas developed independently in \cite{MR666160} by  W.~C.~Hsiang and R.~Staffeldt are tangent to Connes and Tsygan's, but the goal is different: Hsiang and Staffeldt prove that rationally Waldhausen's algebraic K-theory of simply connected spaces can be described by the homology of Lie algebras and the K-theory of $\Z$.  Since we're working rationally, the K-theory of $\Z$ is known from A.~Borel's calculation \cite{MR387496}.


Another indication that things were under control over the rationals came in the form of Goodwillie's
confirmation \cite{MR786354} of Waldhausen's conjecture that stable K-theory and Hochschild homology are rationally equivalent.
This was followed up by \cite{MR855300} showing that
if $A\to B$ is a map of simplicial rings such that $\pi_0A\to\pi_0B$ is a nilpotent extension (\ie a surjective ring homomorphism with nilpotent kernel), then the \emph{relative} K-theory $K(A\to B)=\text{fiber}\{K(A)\to K(B)\}$ agrees with the (shifted) relative cyclic homology 
\begin{equation}
  \label{eq:Goorel}
  K(A\to B)\simeq_{\Q}\Sigma\HH(A\to B)_{h\T}.
\end{equation}
%
Goodwillie's result subsumes that of Hsiang and Staffeldt in the sense that they study the rational K-theory of the $1$-connected map of spherical group rings $\ess[\Omega X]\to\ess$ for $X$ a simply connected space.  However, since K-theory commutes with rationalization this amounts to studying the map $\Q[\Omega X]\to\Q$ which fits in Goodwillie's framework (phrased in terms of simplicial rings after strictifying the loop group into a simplicial group \eg via the Kan loop group).

A crucial ingredient in Goodwillie's analysis is that the Dennis trace factors over the homotopy \emph{fixed} points $\HH(A)^{h\T}$ of the circle action on Hochschild homology
$$K(A)\to\HH(A)^{h\T}\to\HH(A),$$
(see \cite{MR793184} and C.~H.~Hood and D.~S.~D.~Jones \cite{MR920950}) and the equivalence \eqref{eq:Goorel} is induced by this factorization and by the so-called ``norm map'' $\Sigma\HH(A)_{h\T}\to\HH(A)^{h\T}$.  The homotopy groups of the homotopy fixed points are often referred to as ``negative cyclic homology''. 
The positive/negative point of view comes from the grading on the bicomplexes used to calculate (negative) cyclic homology, mirroring the skeleton filtration of the classifying space $B\T=\C P^\infty$ of the circle.

The cofiber of the norm map $\Sigma\HH(A)_{h\T}\to\HH(A)^{h\T}$ (denoted  $\HH(A)^{t\T}$ and called the \emph{Tate construction} or \emph{periodic cyclic homology}) plays an important r\^ole in many situations, and its siblings will reappear below.  For now it suffices to notice that Goodwillie proved that evaluating at $t=0$ gives an equivalence $\HH(A[t])^{t\T}\we\HH(A)^{t\T}$ for rational $A$.

With this development, calculations of rational algebraic K-theory all of a sudden seem feasible, spawning an interest in relevant calculations of cyclic homology; a somewhat arbitrary list could contain papers by D.~Burghelea, G.~Corti\~{n}as, Z.~Fiedorowicz, S.~Geller, L.~Reid, M.~Vigu{\'e}-Poirrier, Villamayor and Weibel \cite{MR834279}, \cite{MR842427},  
\cite{MR1101119}, \cite{MR972360} and certainly a lot of others, some of which will resurface later.
\subsection{Interpretations of the trace map}
\label{sec:Htinvariance}
Connes and Karoubi \cite{ConnesKaroubi} and Weibel \cite{MR902784} provide a good way of expressing the relation between various forms of K-theory and of cyclic homology.  In motivic theory one of the important axioms is $\mathbb A^1$-invariance: the affine line should be contractible.  In the affine situation this means that the theory shouldn't see the difference between a ring $A$ and its polynomial ring $A[t]$.
However, unless you assume that $A$ is particularly nice (for instance regular), $K(A[t])$ and $K(A)$ will be very different, and the difference is given by the (relative) K-theory of the category of nilpotent endomorphisms of finitely generated projective modules.

There is a universal way of imposing $\mathbb A^1$-invariance on a functor $F$ from simplicial rings to spectra: let $HF(A)$ be what you get if you apply $F$ degreewise to the simplicial ring $\Delta A=\{[q]\mapsto A[t_0,t_1,\dots,t_q]/\sum t_j=1\}$.  Then $HF(A[t])\we HF(A)$ and in the case where already $F(A[t])\we F(A)$ the canonical map $F(A)\to HF(A)$ is an equivalence.  Applying this construction to K-theory you get $\mathbb A^1$-invariant K-theory with history going back to Karoubi and O.~Villamayor \cite{MR0313360} and we may define nil-K-theory $nilK(A)$ to be the fiber of the canonical map $K(A)\to HK(A)$. When the degreewise extension of $F$ is a homotopy functor, the contractibility of $\Delta A$ forces $HF(A)$ to be contractible too; this is for instance the case with cyclic homology $\HH(A)_{h\T}$.  Hence, for $A$ rational  the Tate-theory $\HH(A)^{t\T}$ is the $\mathbb A^1$ invariant theory associated with the homotopy fixed points $\HH(A)^{h\T}$ and
the trace gives a map of fiber sequences
$$\xymatrix{nilK(A)\ar[r]\ar[d]&K(A)\ar[r]\ar[d]&HK(A)\ar[d]\\
  \Sigma\HH(A)_{h\T}\ar[r]&\HH(A)^{h\T}\ar[r]&\HH(A)^{t\T}.
}$$
In the relative nilpotent situation of Goodwillie's theorem \eqref{eq:Goorel} the spectra to the right both vanish and in the absolute situation the perspective of the trace as the logarithm from a Lie group to its tangent space (in the guise of cyclic homology) propagandized for above should be viewed as a statement about the left vertical map.


\section{`Rationally' is \emph{not} enough: $\THH$, $\TC$}
\label{sec:BHM}

While all of this was happening, M.~B\"okstedt managed to give a convincing definition of topological Hochschild homology, $\THH$, conforming with Goodwillie and Waldhausen's expectations of a Hochschild homology over the sphere spectrum.
Even though his definition predates the modern view on ring spectra, what B\"okstedt does is to write down an explicit description
(in terms of objects later recognized as good building blocks for (highly structured) ring spectra) closely mimicking the algebraic construction of Hochschild homology  and show that the result has the correct homotopical properties.
Perhaps most impressively, he was even able to calculate the pivotal cases $\THH(\F_p)$ and $\THH(\Z)$.  B\"okstedt's calculations showed that the homotopy ring $\THH_*(\F_p)$ is polynomial on a generator in degree two, starkly contrasting with (derived) Hochschild homology $\HH(\F_p)$ (a.~k.~a. Shukla homology) which is a divided power algebra on a generator in degree two.  So, although $\THH(\F_p)$ and $\HH(\F_p)$ have the same homotopy \emph{groups}, the \emph{ring structures} are vastly different, a phenomenon that is at the heart of many of the ensuing calculations of K-theory through trace methods.

M.~Jibladze, T.~Pirashvili and Waldhausen soon realized that $\THH$ of ordinary rings coincides with an old homology theory of Mac~Lane's and also with the homology of the category of projective modules and of certain functor categories \cite{MR1094244}, \cite{MR1181095}.
This correspondence led to some extremely efficient calculations by V.~Franjou, J.~Lannes, L.~Schwartz and Pirashvili \cite{MR1262942}, \cite{MR1480880}.  In hindsight, some of these groups were already known to L.~Breen \cite{MR516914} through similar considerations.
\subsection{The algebraic K-theory Novikov conjecture}
\label{sec:novikov}
The potential of B\"okstedt's construction was immediately demonstrated through B\"okstedt, Hsiang and Madsen's successful attack \cite{MR1202133} in the late 80's  on the algebraic K-theory version of the Novikov conjecture.  They prove that the assembly map
$$K(\Z)\smsh BG_+\to K(\Z[G])$$
to the K-theory of the integral group ring is rationally injective whenever the group $G$ has finitely generated homology in each dimension.  To achieve this they define ``topological cyclic homology'' $\TC$.
Connes discovered that the Hochschild complex is ``cyclic'' and so (its realization) comes with a natural action by the circle $\T$.
The same goes for topological Hochschild homology, but $\TC$ is not just a homotopy orbit/fixed point-type construction applied to $\THH$.  With the explicit construction B\"okstedt gave, there are fairly apparent maps given in terms of restricting equivariant maps to fixed points.
These restriction maps give a structure on the various fixed point spectra akin to Witt-ring constructions in algebra, a feature of which is captured by saying the $\THH$ is a ``cyclotomic spectrum''. The important ushot is a map $\widehat{\Gamma}_1\colon\THH\to\THH^{tC_p}$.  The actual construction of $\TC$ is somewhat involved, and we cheat by citing the much later interpretation \cite{MR3904731}
of T.~Nikolaus and P.~Scholze, where $\TC$ is the homotopy equalizer
$$\xymatrix{\TC(A)\ar[r]&\THH(A)^{h\T}\ar@<1.2ex>[rr]^{\text{can}}\ar@<-1.2ex>[rr]^{\varphi}&&\THH(A)^{t\T}{\:}\widehat{},}
$$
where $\text{can}\colon\THH(A)^{h\T}\to\THH(A)^{t\T}{\:}\widehat{}{\,}$  is the profinite completion of the cofiber of the above mentioned norm map $\Sigma\THH(A)_{h\T}\to\THH(A)^{hT}$
and the Frobenius $\varphi$ comes from the cyclotomic structure
(or more generally, from the fact that smash powers support ``geometric diagonals'').  Note that if we ignore the profinite completion and replace $\varphi$ by the trivial map, $\TC(A)$ is replaced by (a suspension of) the orbit spectrum $\THH(A)_{h\T}$ -- the most na\"\i ve definition imaginable of ``topological cyclic homology''.

The Dennis trace factors over $\TC$ giving rise to the \emph{(cyclotomic) trace}
$$K(A)\to\TC(A)$$
(this is amazing and it wouldn't have worked with the ``na\"\i ve $\TC$'': the composite $K(A)\to\TC(A)\to\THH(A)^{h\T}$ does not factor over $\Sigma\THH(A)_{h\T}$).

The projection $\ess\to\Z$ induces a rational equivalence $K(\ess[G])\we_{\Q}K(\Z[G])$, so for B\"okstedt, Hsiang and Madsen's purposes it is just as good to work over the sphere spectrum, which is advantageous since they understand $\TC(\ess[G])$ in terms of free loop spaces and power maps
(when applied to spherical group rings the above mentioned restriction maps split.
See also the precursor \cite{MR899917} and the later \cite{MR1408537}).  This means that they understand the assembly map for $\TC$
$$\TC(\ess)\smsh BG_+\to\TC(\ess[G])$$
and eventually the K-theory version of the Novikov conjecture is ``reduced'' to understanding the difference between K-theory and $\TC$ of the ($p$-adic) integers.  In hindsight we know (see Section~\ref{sec:calcinTC}) that after profinite completion the homotopy groups of these spectra agree in non-negative dimensions, but in \cite{MR1202133} they get through knowing somewhat less.

\subsection{Stable K-theory is $\THH$}
\label{sec:KS}

Once B\"okstedt had defined $\THH$, Waldhausen proposed a program to prove that stable K-theory and topological Hochschild homology agree for all $\ess$-algebras, see  \cite{MR1409624}.  Indeed, he had already established the ``vanishing of the mystery homology theory'' by geometric means \cite{MR921486}, which amounts to the equivalence $\ess\we K^S(\ess)$.  Together with the obvious $\ess\we \THH(\ess)$ one has that stable K-theory and $\THH$ coincide in the initial case.

The first full proof was almost orthogonal to this program: Dundas and McCarthy proved $K^S\simeq\THH$ for simplicial rings \cite{MR1307900} by comparing both theories to the homology of the category of finitely generated projective modules through McCarthy's setup for amalgamating (topological) Hochschild homology and Waldhausen's K-theory. The equivalence was later extended to all connective $\ess$-algebra through a denseness argument \cite{MR1607556}.

\section{The first calculation of $\TC$}
\label{sec:BMiandii}
With $\TC$ established as a success for obtaining information about K-theory, B\"okstedt and Madsen attack the most prestigious calculation: What is $\TC(\Z)$?
Their (successful) attack was until recently 
the standard model for next to all calculations of $\TC$ (and -- together with motivic methods -- for K-theory).  The papers \cite{MR1317117} and \cite{MR1368652} are hard to read, but a treasure trove of ideas.
Somewhat streamlined, the procedure is as follows:
\begin{framed}
  \begin{enumerate}[itemsep=0pt]
\item Calculate $\THH$ (B\"okstedt had already done that for $\Z$ and $\F_p$)
\item Consider cartesian squares comparing categorical and homotopical fixed points of the cyclic groups $C_{p^n}$ of orders powers of the prime $p$
  $$\xymatrix{
    \THH^{C_{p^n}}\ar[d]^{\Gamma_n}\ar[r]^-R&\THH^{C_{p^{n-1}}}\ar[d]^{\widehat{\Gamma}_n}\\
    \THH^{hC_{p^n}}\ar[r]&\,\THH^{tC_{p^{n}}.}
  }$$
 
\item Prove that $\widehat{\Gamma}_1$ is an equivalence in high degrees.  
\item Calculate the homotopy fixed points $\THH^{hC_{p^n}}$ and $\THH^{h\T}$ keeping track of known elements coming from K-theory using fixed point/Tate spectral sequences
\item Assemble the pieces keeping in mind possible Witt-like extensions to obtain $\TC$.
\end{enumerate}
\end{framed}
B\"okstedt and Madsen do this at odd primes for the integers, and Rognes does $p=2$ \cite{MR1663390}.
These calculations expose periodic phenomena as predicted by the Lichtenbaum-Quillen conjecture and later expanded to a ``red-shift'' paradigm, see Section~\ref{sec:red}.

Item 2 in the procedure above is not phrased like this in the earliest papers. Nikolaus and Scholze use the diagram to ``define'' the categorical fixed points by means of homotopy invariant constructions, using only the existence of a map $\widehat{\Gamma}_1$ which arises from the cyclotomic structure.  As a matter of fact, their repackaging of the theory has led to new ways of interpreting this attack which seem very fruitful, and entirely new and elegant methods for handling $\TC$ have been developed.  For a striking example see \eg R.~Liu and G.~Wang's preprint \cite{https://doi.org/10.48550/arxiv.2012.15014}.

The reason that proving that $\widehat{\Gamma}_1$ is an equivalence in high degrees suffices is Tsalidis' result \cite{MR1428060} saying that then this is also true for all the $\Gamma_n$'s via an argument akin to Carlsson's proof of the Segal conjecture.


\section{How far apart are K-theory and $\TC$?}
\label{sec:gooICM}

Inspired by his calculus of functor and his rational results, Goodwillie conjectured at the 1990 ICM in Kyoto \cite{MR1159249} that:
\begin{framed}
  \begin{quote}
    Let $K^{\Kinv}$ be the fiber of the trace $K\to\TC$.
    If $A\to B$ is a map of connective ring spectra such that $\pi_0A\to\pi_0B$ is a nilpotent extension, then
    $$K^{\Kinv}A\to K^{\Kinv}B$$
    is an equivalence.
\end{quote}
\end{framed}
The conjecture was stated only for $1$-connected maps, and
Goodwillie's vision was that one should prove that
\begin{enumerate}[itemsep=0pt]
\item Both K-theory and $\TC$ are analytic functors
\item The trace induces an equivalence of differentials
  $$D_1K(A)\we D_1\TC(A).$$
\end{enumerate}
 According to Goodwillie's calculus of functors one would then have that
\begin{quote}
  Within the ``radius of convergence'' the fiber $K^{\Kinv}(A)$ of the trace
  $K(A)\to\TC(A)$
  is constant.
\end{quote}

In the guise of stable K-theory we've seen that $D_1K(A)$ is equivalent to topological Hochschild homology and Hesselholt \cite{MR1317119} proves that so is $D_1\TC(A)$.  From this McCarthy \cite{MR1607555} proves Goodwillie's conjecture after profinite completion for simplicial rings, and this was extended to connective $\ess$-algebras in \cite{MR1607556}.  The integral result 
stated above as Goodwillie's conjecture was established around 1996 and documented in \cite{MR3013261}.

\section{Calculations of $\TC$ amount to calculations in K-theory}
\label{sec:calcinTC}

Algebraic K-theory was already fairly well understood ``away from the prime $p$'' by results like Gabber rigidity \cite{MR1156502} which, for instance, implies that the projection $\Z_p\to\F_p$ from the $p$-adics to the prime field induces an equivalence in K-theory after completion at any prime different from $p$.  What topological cyclic homology has to offer is that we get hold of information ``at the prime''.

It all starts with the prime field $\F_p$, which Quillen tells us is fairly dull at $p$:
$$K(\F_p)\simeq_p H\Z,$$
where $\simeq_p$ designates an equivalence after $p$-completion and $H\Z$ is the integral Eilenberg-Mac-Lane spectrum.  Madsen \cite{MR1341845} shows that the cyclotomic trace is an equivalence except for noise in negative degrees: 
$\TC(\F_p)\simeq_p\Sigma^{-1}H\Z\vee H\Z.$
Goodwillie's conjecture then implies that for all $n>0$
$$K(\Z/p^n\Z)\simeq_p\TC(\Z/p^n\Z)[0,\infty),$$
where the spectrum on the right hand side is the connective cover.  Now, investigations by, for instance, A.~Suslin and I.~Panin \cite{MR864175} show that after profinite completion,  $K(\Z_p)$ is the limit of the $K(\Z/p^n\Z)$s, and when Hesselholt and Madsen \cite{MR1410465} establish the same result for $\TC$ they are in position to conclude that B\"okstedt and Madsen's calculation of $\TC(\Z)$ amounts to a calculation of the $p$-completion of $K(\Z_p)$!

On the $\TC$ side, \cite{MR1410465} cements and enhances the innovations of \cite{MR1317117} and \cite{MR1368652}.  The paper also provides a lot of valuable perspectives (plus that the end result is more general than I've stated) and should be seen in tandem with \cite{MR1417085} which realizes Bloch's program connecting K-theory and crystalline cohomology \cite{MR488288}.

At the end of the millennium, the best structural result may be summarized by
\begin{framed}
  \begin{quote}
    Let $K^{\Kinv}$ be the fiber of the trace $K\to\TC$.
    If $k$ is a perfect field of characteristic $p>0$, $A$ a connective ring spectrum whose ring of path components $\pi_0(A)$ is a $W(k)$-algebra which is finitely generated as a module, then
    $$K^{\Kinv}(A)\simeq_p\Sigma^{-2}H\mathrm{coker}\{\pi_0A\overset{1-F}{\longrightarrow}\pi_0A\},$$
  where $F$ is the Frobenius. 
\end{quote}
\end{framed}

\section{Where do we go from here? $\TC$ in the $21$st century.}
\label{sec:wheredowegofromhere}
The pioneering age of topological cyclic homology ends around the turn of the millennium.  In the $21$st century $\TC$ consolidates its r\^ole as a work horse in the study  K-theory.
Firstly, $\TC$ provides calculations of K-theory that seemed totally out of reach in 1990. We give some examples of this in Section~\ref{sec:explcalc}.

Secondly, after 2000 the structural understanding of $\TC$ expands immensely, shedding light on the nature of $\TC$ as well as on K-theory and we give a brief account in Section~\ref{sec:struc}.   Many of the structural results were pivotal to the calculations; for instance, understanding how $\TC$ behaves under localization was needed for calculations of the K-theory of local number fields.  This structural understanding has even spilled over to neighboring fields -- a technical part of the calculation of K-theory of the integers via motivic cohomology as presented by Rognes and Weibel \cite{MR1697095} relied on understanding the multiplicative structure of $\TC$. 

Last but not least, new frameworks and perspectives on $\TC$ and the trace have materialized; a much too brief discussion is concentrated in Section~\ref{sec:new}.



\section{An explosion of calculations}
\label{sec:explcalc}

\subsection{Local fields}
\label{sec:numberlqc}

In \cite{MR1998478} Hesselholt and Madsen muster the full force of the theory to calculate the K-theory of local fields.
To achieve this they use an amalgam of
\begin{itemize}[itemsep=0pt]
\item the methods for calculation of $\TC$ initiated by B\"okstedt, Madsen and Hesselholt, enhanced through setting up a framework for handling the ``de~Rham-Witt complex''
\item Goodwillie's conjecture
\item a setup for handling localization through $\log$-formalism combined with $\TC$ of categories -- a precursor of  $\log$-methods and localization for ring spectra.
\end{itemize}
The last point is crucial and we will return to a fuller discussion of  $\log$-methods and localization for ring spectra in Section~\ref{sec:localization}.  The calculation of the K-theory of local fields \cite{MR1998478} accentuates how structurally different K-theory and $\TC$ are, in particular how differently they behave with respect to localization.  Whereas $\TC$ of a field of characteristic $0$ is pretty uninformative in positive degrees, the K-theory of a local field $F$ sits in a fiber sequence $K(k)\to K(\mathcal O_F)\to K(F)$ where $\mathcal O_F$ is the ring of integers and $k$ is the residue field.  The failure is not so much $\TC$ not respecting Waldhausen's localization sequence as $\TC$ not allowing a reformulation of the base spectrum as $\TC(F)$.   Taking the consequence of this, Hesselholt and Madsen simply keep the category Waldhausen provides and take $\TC$ of that using  \cite{MR1388700} to obtain a map from the localization sequence in K-theory.  See also \cite{MR1626647} for a structural comparison of  K-theory and $\TC$. 

From an algebraic and equivariant point of view, the treatment of the de~Rham-Witt complex may be even more interesting.  In the paper this mixture of the ring of Witt vectors and the de~Rham complex (with $\log$-poles) is important because it provides a framework for handling the interplay between the various fixed point spectra of $\THH$.

The setup has been revisited many times, for instance in Hesselholt, M.~Larsen and A.~Lindenstrauss' \cite{MR4050106}.

A reason the methods work for local but not global fields is exemplified by the observation that $K(\Z)$ and $K(\Z_p)$ are vastly (and interestingly) different, while $\TC(\Z)\we_p\TC(\Z_p)$.

\subsection{Truncated polynomials, nilterms, endomorphisms}
\label{sec:trunc}
Calculations of K-theory of truncated polynomial rings go far back to the infancy of K-theory, and with Goodwillie's conjecture an established fact such calculations became open for attack; the first published account being \cite{MR1471886}  where the answer is given in terms of big Witt vectors.   This topic have been followed up by several authors, for instance V.~Angeltveit's \cite{MR3607214}.  Also, there is an interesting interplay between $\mathbb A^1$-invariance and nil-invariance, making calculations of the K-theory of nilpotent endomorphism somewhat more accessible.

The identification between stable K-theory and $\THH$ factors over the K-theory of the category of endomorphisms, which for a very long time had been known to be closely related to the theory of Witt vectors.  Before the existence of $\TC$, Goodwillie had written a widely circulated letter to Waldhausen indicating a theory in this direction.  McCarthy soon realized that there was a close connection between the Taylor tower of K-theory and such a Witt-like theory and saw it as an extension of the $\TR$-construction used to build $\TC$, see \eg Lindenstrauss, McCarthy \cite{MR2928981},    {A.~Blumberg, D.~Gepner, G.~Tabuada \cite{MR3430369} and the recent \cite{MR4381922} by C.~Malkiewich and K.~Ponto.}

The vanishing of negative K-groups is a somewhat related question that has been addressed through trace methods.   
For the localization square
$$\xymatrix{K(\mathbb P_A)\ar[r]\ar[d]&K(A[t])\ar[d]\\K(A[t^{-1}])\ar[r]&K(A[t,t^{-1}])}$$
to be cartesian (where $\mathbb P_A$ is the projective line in Quillen's sense; $K(\mathbb P_A)\simeq K(A)\vee K(A)$ \cite{MR0338129}), Bass observed that you need to extend the definition of K-theory to also have negative homotopy groups (see \eg \cite{MR1106918}) and Weibel formulated a question about whether the negative K-groups vanish below the dimension of the input.  Using (topological) cyclic homology this was shown in characteristic zero by in \cite{MR2415380} by G.~Corti\~nas, C.~Haesemeyer, M.~Schlichting and Weibel
and in positive characteristic (over infinite perfect fields, assuming strong resolution of singularities) by T.~Geisser and Hesselholt in \cite{MR2677901}.

\subsection{Algebraic K-theory of spaces, ring spectra and assembly}
\label{sec:grouprings}
Before Grothendieck, J.~H.~C.~Whitehead developed K-theoretical methods for group rings, and to this day (spherical) group rings are among the most important inputs for K-theory via  Waldhausen's algebraic K-theory of spaces as discussed briefly in Section~\ref{sec:prehistory}.   Goodwillie's conjecture opened a new line of investigation.  For an overview of the situation and some more background,  see Rognes' very readable lecture \cite{MR1743242} where he also outlines the link between automorphisms of manifolds and algebraic K-theory in more detail.

As an example let's look at the initial case: $\ess\to\Z$.  As promised, $K(\ess)$ contains much geometric information: indeed, according to Waldhausen there is an isomorphism
$$K_i(\mathbb S)\cong \pi_i\mathbb S\oplus \mathrm{colim}_n\pi_{i-2}\mathrm{Diff}(D^{n+1}\,\mathrm{ rel }\, D^n),$$
where $\mathrm{Diff}(D^{n+1}\,\mathrm{ rel }\, D^n)$ is the space of diffeomorphisms of the closed unit $(n+1)$-dimensional disc that are fixed on the Northern hemisphere.  The colimit is actually attained in the sense that each homotopy class of such diffeomorphisms of discs of sufficiently high dimension is present in the algebraic K-theory of the sphere spectrum.  On the other hand, K-theory detects a lot of number theory, starting with class field theory in dimension zero and continuing in all dimensions as codified for instance in the Lichtenbaum--Quillen conjecture. 

So, the induced map $K(\ess)\to K(\Z)$ spans all the way from information about diffeomorphisms of 
discs to number theory, but the chasm is not wider than what $\TC$ can handle.   Goodwillie's conjecture tells us that the square
$$\xymatrix{
  K(\ess)\ar[r]\ar[d]&\TC(\ess)\ar[d]
\\
K(\Z)\ar[r]&\TC(\Z)
}$$
is cartesian.  Add to that the calculation of $\TC(\Z)$ and the homotopy theoretic control of $\TC(\ess)$
we see that we're in an extremely interesting situation.
%
J.~Klein and Rognes attacked the problem 
in \cite{MR1432423} and the issue has been revisited by many since,
see \eg Rognes' identification \cite{MR1988283} of the homotopy type of $K(\ess)$ 
and Blumberg and M.~Mandell's complement on the homotopy groups in \cite{MR3921317}.

The next space one wants to understand the algebraic K-theory of is the circle $S^1$.
As is commented by Madsen in \cite{MR1341845},  thanks to the insight provided by Waldhausen, K.~Igusa \cite{MR972368}, F.~T.~Farrell and L.~E.~Jones \cite{MR1159252} and M.~Weiss and R.~Williams \cite{MR953917}, knowing the cofiber of the assembly map $K(\ess)\smsh S^1_+\to K(\ess[\Omega S^1])=K(\ess[t,t^{-1}])$ (which is called the topological Whitehead spectrum $Wh^\mathrm{TOP}(S^1)$ of the circle) goes a long way towards understanding the homeomorphisms of manifolds with negative sectional curvature.  For further details, see \eg the introduction to \cite{MR2597738} where Hesselholt provides some calculations which should be seen in conjunction with J.~Grunewald, Klein and T.~Macko's \cite{MR2399133}.


As we saw, $\TC$ was created to solve the Novikov conjecture.  There is a host of refined conjectures and results in this direction -- essentially about how $K(A[G])$ can be decomposed into information about $K(A)$ and information about the subgroup lattice of $G$ -- but just some of the progress can be attributed $\TC$.  This progress can be followed by consulting the paper \cite{MR4015233} by W.~L\"uck, H.~Reich, Rognes and M.~Varisco.

Ausoni and Rognes' calculations of the algebraic K-theory of topological K-theory \cite{MR1947457} (to be commented further on in the next section) is another prime example of calculations of K-theory of ring spectra through trace methods.  Later this calculation was given a geometric interpretation through two-vector bundles \cite{MR2832571}, \cite{MR3010546}, but the concrete calculation needed to show that the Dirac monopole splits as a virtual two vector bundle \cite{MR2466184} was only accessible through a comparison using the trace.

\label{sec:thom}
Many important ring spectra are Thom spectra.  Blumberg, C.~Schlichtkrull and R.~Cohen \cite{MR2651551} show that Thom spectra are very manageable on the $\THH$-side -- a far-reaching extension of the equivalence between $\THH(\ess[G])$ and the free loop space $\mathrm{Map}(S^1,BG)$.  Most provocatively, $H\F_p$ is a Thom spectrum (but admittedly not through an $E_\infty$-map). Following the references to \cite{MR2651551} will lead to many recent results where familiar spectra are viewed as Thom spectra, and, among other things, to reformulations of B\"okstedt's calculation, see \eg the multiplicative version in A.~Krause and Nikolaus \cite{KN}.

Another amazing equivariant result about $\THH$ of Thom spectra comes in the form of S. Lun{\o}e-Nielsen and Rognes' ``Segal conjecture'' for complex cobordism $\mathrm{MU}$ \cite{MR2832570}.  This makes it attractive to try to understand algebraic K-theory of spaces not via descent over $\ess\to\Z$, but over $\ess\to\mathrm{MU}$, see Section~\ref{sec:localization}.

See also Malkiewich \cite{MR3686399} for a dual picture.

\subsection{Red-shift/chromatic behavior beyond LQC}
\label{sec:red}

Chromatic convergence (from now on: complete everything at a prime $p$ and ignore fine points wrt. what localization to use) tells us that the sphere spectrum is the limit of the tower
$$\dots\to L_n\ess\to L_{n-1}\ess\to\dots\to L_1\ess\to L_0\ess=H\Q$$
of localizations at the Johnson-Wilson theories $E(n)$ with coefficients $E(n)_*=\Z_{(p)}[v_1,\dots,v_n,v_n^{-1}]$.  Waldhausen \cite{MR764579} speculated that applying K-theory would yield a tower starting at $K(\Q)$ and converging to $K(\ess)$ -- giving a highly structured ascent from number theory to manifolds!  See J.~McClure and Staffeldt \cite{MR1164148} for more details.

C.~Ausoni and Rognes take the first step up this ladder beyond the Lichtenbaum-Quillen conjecture in  \cite{MR1947457} where they calculate the K-theory of the so-called Adams summand; followed up by Ausoni's calculation \cite{MR2609252} of the K-theory of the connective complex K-theory, $\mathrm{ku}$.  Part of the pattern that emerges can be seen from the following table
(think of $\Z_p$ as ``$\F_p[[v_0]]$''):
\begin{table}[htbp]
  \centering
  \begin{tabular}{|r|c|l|}
    \hline
    the homotopy of &$\TC(\F_p)$&is a finitely generated free $\Z_p$-module\\
    the mod-$p$  homotopy of &$\TC(\Z)$&is a finitely generated free $\F_p[v_1]$-module\\
    the mod-$(p,v_1)$  homotopy of &$\TC(\mathrm{ku})$&is a finitely generated free $\F_p[v_2]$-module\\
                                                  \hline
  \end{tabular}
\end{table}
\newline\noindent
Ausoni and Rognes also point to conjectural calculations showing a clear pattern, but depending on properties of the Brown-Petersen spectra that must be somewhat adjusted.  They synthesize the picture into a ``red-shift conjecture'' in \cite{MR2499538}.  See also \cite{MR3728661} which among many other things points towards a motivic filtration for commutative ring spectra, see \eg \cite{https://doi.org/10.48550/arxiv.2206.11208}.

The chromatic properties of K-theory has been revisited several times, \eg B.~Bhatt, D.~Clausen, A.~Mathew \cite{MR4110725}, Ausoni, B.~Richter and Rognes \cite{MR4071375}, \cite{MR2928844}. Personally, my favorite approach to red-shift is through iterations, noting that (for a commutative ring spectrum) iterating $\THH$ is simply tensoring with a torus whose rich group of symmetries should structure the chromatic behavior.  A first start can be found in T.~Veen's thesis \cite{MR3748678} which shows that the periodic self-map $v_{n-1}$ can be detected in the $n$-fold iterated K-theory $K^{(n)}(\F_p)$, at least for $n\leq p>3$.

Crucial for this way of thinking is the extension of the ``cyclotomic'' concept for $\THH$ alluded to in Section~\ref{sec:novikov} to the realization that the geometric fixed points of $X\otimes A$, where $X$ is a $G$-space and $A$ a commutative ring spectrum, is given by the homotopy $G$-orbits $X_{hG}\otimes A$ through the ``geometric diagonal'' \cite{MR2729005}, \cite{BDS} ($X\otimes A$ is the categorical tensor in commutative ring spectra -- a notion that already appears in J.~McClure, R.~Schw{\"a}nzl and R.~Vogt's interpretation of topological Hochschild homology \cite{MR1473888}).
Put another way, we really don't have a good picture for why there should be red-shift for K-theory, but the very tight connection between K-theory and $\TC$ offers us the alternative -- and to me reasonable -- explanation that it is is really about the symmetries in stable homotopy theory to which K-theory is just an innocent victim.
\footnote{{\bf Stop the press!} Since this talk was held, the history of the red-shift conjecture has been seriously altered and most likely even come to a satisfactory conclusion for commutative ring spectra with the preprint on the chromatic Nullstellensatz by R.~Burklund, T.~Schlank and A.~Yuan \cite{https://doi.org/10.48550/arxiv.2207.09929}.  Precisely, their claim is (Theorem~9.11)
  \begin{quote}
    Let $A$ be a nontrivial commutative ring spectrum of nonnegative height.  Then the height of $K(A)$ is exactly one greater than the height of $A$. 
  \end{quote}
  This statement relies on J.~Hahn's result \cite{https://doi.org/10.48550/arxiv.1612.04386} on the chromatic support of commutative ring spectra to give a good characterization of the height.  Burklund, Schlank and Yuan need to show that the height of $K(A)$ is strictly greater that the height of $A$ since the fact that  the difference between the heights of $K(A)$ and $A$ is not more than one follows from M.~Land, Mathew, L.~Meier, G.~Tamme \cite{https://doi.org/10.48550/arxiv.2001.10425} and Clausen, Mathew, Nauman, Noel \cite{https://doi.org/10.48550/arxiv.2011.08233}. They achieve this through their Nullstellensatz which guarantees a map to a Lubin-Tate theory where Yuan [Yua21] already has shown the result.

  Note that full commutativity is essential to this setup; what happens in the $E_n$-case for $n<\infty$
  is not understood. Forerunners to this phenomenon were noticed already by Ausoni and Rognes in their conjectural calculations for $BP$-theories in the introduction to \cite{MR1947457}. For recent development in this direction the reader may enjoy Hahn and D.~Wilson's preprint \cite{https://doi.org/10.48550/arxiv.2012.00864}.

  When the increase in height is established, it may be the time to look back at the (admittedly rather weak) calculational evidence that lead to the conjecture and see whether there is something more detailed to be said about the structure of algebraic K-theory (or perhaps rather: about the interplay between group actions and commutativity).


}

\section{Structural results}
\label{sec:struc}

\subsection{Structure {\bf on} $\TC$ 
}
\label{sec:structureON}

Although ``we'' are awesome at calculating $\TC$ we know very little about the deeper structure of $\TC$. For instance, even though $\TC$ preserves much structure, we don't understand the commutative structure of $\TC(A)$ -- or even $\THH(A)$ -- of a commutative ring spectrum $A$.  My favourite example is the following: the commutator $S^1\to S^1\vee S^1$ (attaching the $2$-cell in the torus) induces a map of augmented commutative $H\F_p$-algebras
$$\THH(\F_p)\to\THH(\F_p)\smsh_{H\F_p}\THH(\F_p).$$
Is it trivial?  On homotopy it is the trivial map $\F_p[x]\to\F_p[x_1,x_2]$, $x\mapsto 0$ and furthermore as an \emph{associative} 
map it is trivial  and all primary $E_\infty$-obstructions to triviality vanish.  If it \emph{is} trivial as an $E_\infty$-map, we have a slew of results relevant to red-shift.

That said, we know some things, for instance from H.~Bergsaker and Rognes \cite{MR2602850} and also Angeltveit and Rognes \cite{MR2171809}.  Another point of view relevant for the multiplicative properties is the ``norm'' or ``smash power'' point of view exemplified by \cite{MR2737802} and \cite{MR3933034}  and which we touched upon in Section~\ref{sec:red}.  

\subsection{New frameworks}
\label{sec:new}

We've already mentioned Nikolaus and Scholze's new take on $\TC$.  This perspective is so useful and has proven to be so influential that I've allowed myself to anachronistically refer to it several times already -- even when introducing $\TC$ in Section~\ref{sec:BHM}.
See Hesselholt and Nikolaus' survey \cite{MR4197995} for examples of this framework used on old and new problems.

There is another very important point of view that puts things in a different perspective.  From the early days the universal properties of K-theory was quite apparent: roughly stated as ``algebraic K-theory is the universal theory under the nerve that satisfies Waldhausen's additivity theorem''.  Things were finally firmly organized around this principle with Tabuada's additive and localizing invariants, see \cite{MR3194756} and together with Blumberg and Gepner \cite{MR3070515}, \cite{MR3209352} and also by C.~Barwick \cite{MR3465850}.  In general, the $\infty$-categorical perspective has permeated much of the more recent literature.  An example where the Goodwillie conjecture is significantly extended can be found in E.~Elden's \cite{https://doi.org/10.48550/arxiv.2010.09155}.

\subsection{Back to algebraic geometry}
\label{sec:nonaffine}
Grothendieck designed $K_0$ to answer questions in algebraic geometry and the trace methods were from the start built around Chern characters.  As we've seen, $\TC$ doesn't always behave well with respect to localizations, but in the smooth case many things are fine, and Geisser and Hesselholt set up a successful theory for $\TC$ of schemes in \cite{MR1743237} via hypercohomology as an alternative to $\TC$ of the exact category of algebraic vector bundles.  Along the way they also provide a lot of structural insight and important results.
In subsequent papers Geisser and Hesselholt follow up with a string of impressive calculation, much, but not all, in positive characteristic.  This line of investigation continues to bear fruit, with calculations in algebraic geometry and feeding back to the theory, see for instance \cite{MR2171226}.

Some more recent results in this direction are connected with $p$-adic Hodge theory, see \eg Bhatt, Morrow and Scholze's \cite{MR3905467} and  \cite{MR3949030}.  Scholze's theory of perfectoid rings plays a central r\^ole here, as does Hesselholt's much earlier paper \cite{MR2222509}.  A commutative $\Z_p$-algebra $R$ is perfectoid 
if there exists a nonzero-divisor $\pi\in R$ with  $p\in\pi^pR$ such that $R$ is complete and separated in the $\pi$-adic topology and such
that the Frobenius $\phi\colon R/\pi\to R/\pi$ is a bijection.

As an appetizer, let's mention that Bhatt, Morrow and Scholze extend B\"okstedt periodicity to all perfectoid rings $R$: $\pi_*\THH(R)_p\cong R[u]$ and identify Nikolaus and Scholze's maps $\phi, \mathrm{can}\colon\THH(R)^{h\T}_p\to\THH(R)^{t\T}_p$ on homotopy groups as
$$\xymatrix{
  A_\text{inf}(R)[u,v]/(uv-\xi)
  \ar[rrrr]^-{\phi(u)=\sigma,\,\phi(v)=\phi(\xi)\sigma^{-1}}_-{\text{can}(u)=\xi\sigma,\,\text{can}(v)=\sigma^{-1}}&&&&
  A_\text{inf}(R)[\sigma,\sigma^{-1}].
  }
  $$
  Here $|u|=-|v|=2$ and $A_\text{inf}(R)$ is Fontaine's universal
$p$-complete pro-infinitesimal thickening given as the ring of Witt vectors of the limit of $R/p$ along the Frobenius map and $\xi$ is a generator of the kernel of the associated map $\theta\colon A_\text{inf}(R)\to R$.  A fascinating point in their calculation is where they work not over $\ess$ but the polynomial ring $\ess[z]$ showing that the cyclotomic structure fits equally well here.

This is a point of departure where the interest in the homological side becomes clear, even at the expense of the connections to K-theory, a theme that we unfortunately don't have the occasion to explore here and where I hope someone can write an accessible survey.  Even so, there is payback to K-theory, for instance in the form of S.~Kelly and Morrow's proof of Gersten injectivity for valuation rings over fields \cite{MR4264079}.


\subsection{Localization and Galois theory for ring spectra}
\label{sec:localization}
A central question in K-theory is the following: if a group $G$ acts on a ring $A$ through ring homomorphisms (of particular interest is the case of a Galois extension $A^G\subseteq A$), to what extent is the natural map from the K-theory of the fixed ring to the homotopy fixed points of the K-theory
$$K(A^G)\to K(A)^{hG}$$
an equivalence?  At the beginning of the millennium Rognes embarked on the quest to extend Galois theories to ring spectra (part of the motivation can be found in the introduction of \cite{MR1947457}) and much of this program is documented in \cite{MR2387923}.  This tour de force has seen applications in many direction and is a guiding principle also for the trace invariant.

Ever since the 1980s (see \eg R.~Thomason's perspective on K-theory \cite{MR1106918}) descent has been a central tool.  K-theory satisfies Nisnevich descent, but  not necessarily \'etale descent-- unless one localizes.  Clausen, Mathew, Nauman and Noel address and partially answer the question about to what extent a Galois action passes through K-theory in \cite{MR4071324}.   Successful applications that use both trace methods and descent have often proceeded by proving descent in parallel in K-theory and $\TC$, see for instance  \cite{MR3505658} and \cite{MR4071324}.

Indeed, many applications of trace methods to descent have taken the following form: assume given a cosimplicial ring spectrum $A_\bullet$ (or some other diagram), such that for each $[n]\in\Delta$ you have control over the trace $K(A_{[n]})\to\TC(A_{[n]})$, to what extent does that give you control over the trace on $A=\mathrm{holim}_\Delta A_\bullet$?  In \cite{MR1607556} this was used to transport theorems from simplicial rings to connective ring spectra by considering the Amitzur complex $[n]\mapsto H\Z^{\smsh n+1}$ converging to $\ess$.  While of theoretical importance, the spectra $H\Z^{\smsh n+1}$ are computationally useless, and in \cite[4]{MR3873116} the prospect of using descent along $\ess\to MU$ is suggested.  This has the advantage that $MU^{\smsh n+1}$ is more accessible than $H\Z^{\smsh n+1}$ since $MU$ is a nice Thom spectrum , $MU^{\smsh n+1}\simeq MU\smsh BU^{\times n}_+$ and the Segal-type control of the equivariant structure offered by \cite{MR2832570}.

The reasons that prompted the combination of trace and of $\log$-methods in \cite{MR1998478} continue to be relevant when considering ring spectra. This need was first addressed by Rognes in \cite{MR2544395} and has since been refined and followed up, see \eg Rognes, S.~Sagave and Schlichtkrull \cite{MR2964635}, \cite{MR3412362}.

Another need that became clear with \cite{MR1947457} was a good theory for localizations for ring spectra, and the first test was whether there is a map of localization sequences
$$\xymatrix{K(\Z)\ar[r]\ar[d]& K(ku)\ar[r]\ar[d]& K(KU)\ar[d]\\\TC(\Z)\ar[r]& \TC(ku)\ar[r]& \TC(ku\mid KU),}
$$
where $KU$ is the localization of $ku$ you get by inverting the Bott element $u\in ku_2$.  This challenge was taken up by Blumberg and Mandell who establish this in \cite{MR2413133} and \cite{MR2928988}.  Unfortunately, this picture does not extend beyond $ku$ to higher Brown-Peterson spectra, see  B.~Antieau, T.~Barthel and Gepner \cite{MR3760300} who use the trace to show that the fiber of localization can't be as hoped.  That said, according to \cite[Theorem 0.3]{MR3760300}, if $A$ is a ring spectrum and $a\in\pi_tA$ is an element such that $\{1,a,a^2,\dots\}$ satisfies the Ore~condition, then the localization $A\to A[a^{-1}]$ gives rise to a fiber sequence
$$K(End_A(A/a)^\text{op})\to K(A)\to K(A[a^{-1}]).$$

\subsection{Singularities/closed excision}
\label{sec:excision}
Except for the applications of $\TC$ to nilpotent extensions, the initial development focused on the smooth case (nicely complemented by motivic methods).  However, for K-theory ``closed excision'' was important from the start.  In the affine case it can be stated as follows: consider a ``Milnor square'', a cartesian square of rings 
$$\xymatrix{A\ar[r]\ar[d]&B\ar@{->>}[d]\\C\ar[r]&\,D;}$$
where $B\fib D$ is surjective.  What can you say about $K(A)$?  In the rational case this received quite some attention in the eighties, but it was not before Corti\~nas \cite{MR2207785}, Geisser and Hesselholt \cite{MR2249803}, Dundas and H.~Kittang \cite{MR2369021}, \cite{MR3031812} and M.~Land and G.~Tamme \cite{MR4024564} that the full result was established
\begin{framed}
  \begin{quote}
    Let $K^{\Kinv}$ be the fiber of the trace $K\to\TC$.
    If $$\xymatrix{A\ar[r]\ar[d]&B\ar[d]\\C\ar[r]&D}$$ is a cartesian square of connective ring spectra where $\pi_0B\fib\pi_0 D$ is surjective, then the associated square $$\xymatrix{K^{\Kinv}A\ar[r]\ar[d]&K^{\Kinv}B\ar[d]\\K^{\Kinv}C\ar[r]&K^{\Kinv}D}$$
    is also cartesian.
  \end{quote}
\end{framed}
Note that Land and Tamme prove more than closed excision.

As with Goodwillie's conjecture, closed excision opens up the possibility for calculations of K-theory, see for instance Hesselholt, Angeltveit T.~Gerhardt and Nikolaus \cite{MR2301459},   \cite{MR2846160} and \cite{MR3290091}.

\subsection{Rigidity}
\label{sec:rigidity}

I end the talk with a personal note and a cool result with an amazing proof.
In 1996 I applied for a grant to show ``rigidity'' in the following form: If $A\to B$ is a map of connective ring spectra so that $\pi_0A\to \pi_0B$ is a complete extension (a surjection with a complete kernel), then
$$K^{\Kinv}A\to K^{\Kinv}B$$
is an equivalence after profinite completion.
I didn't get the grant (which is my only excuse for not delivering -- though I shared some beers with Randy McCarthy the one evening we both believed that we knew a proof).

In the case when $\pi_0A$ and $\pi_0B$ are commutative there is a generalization of the completeness assumption.  A surjection of commutative rings $f\colon R\to S$ is a henselian extension if for all polynomials $q\in R[x]$ with reduction $q_S\in S[x]$ the induced function
\[\{
  \text{roots $\alpha$ of $q$ such that $q'(\alpha)$ is a unit}\}
  \to \{\text{roots ${\beta}$ of $q_S$ such that ${q_S'}(\beta)$ is a unit}\}
  \]
  is surjective.  
  The proof of rigidity by Clausen, Mathew and Morrow \cite{MR4280864} is truly wonderful.  A crucial ingredient is the fact (which I never would have guessed) that $\TC/p$ commutes with filtered colimits and also that
  the trace mod $p^n$ is split injective on homotopy groups when evaluated at local $\F_p$-algebras.
  \begin{framed}
    \begin{quote}
      Let $K^{\Kinv}$ be the fiber of the trace $K\to\TC$.
      If $A\to B$ is a map of connective ring spectra such that $\pi_0A\to\pi_0B$ is a henselian extension, then
      $$K^{\Kinv}A\to K^{\Kinv}B$$
      is an equivalence after profinite completion.
    \end{quote}
  \end{framed}

  {\bf Exercise.}
    Rigidity shouldn't depend on commutativity assumptions.
    State and prove the correct theorem.

    Redshift should depend on some commutativity.  Explore the borderline cases and explain the fine-structure.

    Merging closed excision with Goodwillie's conjecture allows for more general statements.  Find and prove a satisfactory statement for higher dimensional cartesian cubes.



\bibliographystyle{amsalpha}
\providecommand{\bysame}{\leavevmode\hbox to3em{\hrulefill}\thinspace}
\providecommand{\MR}{\relax\ifhmode\unskip\space\fi MR }
\providecommand{\MRhref}[2]{%
  \href{http://www.ams.org/mathscinet-getitem?mr=#1}{#2}
}
\providecommand{\href}[2]{#2}
\bibliographystyle{amsalpha}
\bibliography{fields.bib}
\end{document}